\documentclass{article}%
\usepackage{amsmath}
\usepackage{amsfonts}
\usepackage{amssymb}
\usepackage{graphicx}%
\providecommand{\U}[1]{\protect \rule{.1in}{.1in}}
\newtheorem{theorem}{Theorem}[section]

\newtheorem{corollary}[theorem]{Corollary}

\newtheorem{definition}[theorem]{Definition}
\newtheorem{example}[theorem]{Example}

\newtheorem{lemma}[theorem]{Lemma}

\newtheorem{proposition}[theorem]{Proposition}
\newtheorem{remark}[theorem]{Remark}

\newenvironment{proof}[1][Proof]{\noindent \textbf{#1.} }{\  \rule{0.5em}{0.5em}}
\numberwithin{equation}{section}
\begin{document}

\title{Modelling Derivatives Pricing Mechanisms with Their Generating Functions}
\author{Shige PENG\thanks{ This research is supported in part by National Natural
Science Foundation of China No. 10131040. and Research Foundation and Shanghai
Future Exchange. \ }\\Institute of Mathematics, Institute of Finance\\Shandong University\\250100, Jinan, China\\School of Mathematics, Fudan University \\peng@sdu.edu.cn}
\date{Version: May 2006}
\maketitle

\bigskip \medskip \medskip

\noindent \textbf{Abstract. }{\small In this paper we study dynamic pricing
mechanisms of financial derivatives. A typical model of such pricing mechanism
is the so-called }$g${\small --expectation} {\small defined by solutions of a
backward stochastic differential equation with }$g$ {\small as its generating
function. Black-Scholes pricing model is a special linear case of this pricing
mechanism. We are mainly concerned with two types of pricing mechanisms in an
option market: the market pricing mechanism through which the market prices of
options are produced, and the ask-bid pricing mechanism operated through the
system of market makers. The later one is a typical nonlinear \ pricing
mechanism. Data of prices produced by these two pricing mechanisms are usually
quoted in an option market. }

{\small We introduce a criteria, i.e., the domination condition (A5) in
(\ref{domination}) to test if a dynamic pricing mechanism under investigation
is a $g$--pricing mechanism. This domination condition was statistically
tested using CME data documents. The result of test is significantly positive.
We also provide some useful characterizations of a pricing mechanism by its
generating function. }

\medskip

\medskip \newpage

\noindent \newline

\section{Introduction\label{ss1}\bigskip}

How to quantitatively describe the pricing mechanism of a market of
derivatives is a very interesting problem. A model of dynamic pricing
mechanism of derivatives is formulated (see (A1)--(A4) in the next section) to
characterize this pricing behavior.

We are mainly concerned with two types of pricing mechanisms in an option
market: the market pricing mechanism which outputs the trading prices of
options and the bid--ask pricing mechanism operated according the system of
market makers. \ We stress here that, in our point of view, the ask prices and
the bid prices quoted in a market are determined by a single pricing
mechanism. The difference of a ask price and the corresponding bid price,
called bid--ask spread, reflects the nonlinearity of this mechanism. The data
of prices of above mentioned two pricing systems is usually systematically
quoted in the internet thus the models under our investigation can be
statistically tested. We hope that our modelling can also be applied to
describe the pricing mechanism of some other financial institutions.

The well-known Black--Scholes formula is a typical model of dynamic pricing
mechanism of derivatives. It is a linear pricing mechanism. In fact, the
prices produced by this mechanism is solved by a linear Backward Stochastic
Differential Equation (BSDE). This means that the corresponding generating
function $g$ of the BSDE is a linear function. Nonlinear pricing model by BSDE
was proposed in \cite{EPQ1997} (cf. \cite{EK-Q}). In this paper we show that
each well-defined BSDE with a fixed generating function $g$ forms a dynamic
pricing mechanism, called $g$--expectation and that the behaviors of this
mechanism are perfectly characterized by the behaviors of $g$. Several
conditions of equivalence provided in this paper will be very helpful to
characterize and to find the generating function, or in some other
circumstances, to regulate or to design a pricing mechanism.

A very interesting problem is how to design a test procedure to verify whether
an existing pricing mechanism of derivatives is a $g$--expectation. We will
present the following result: if a dynamic pricing mechanism is uniformly
dominated by a $g_{\mu}$--expectation with a sufficiently large number $\mu$
for the function $g_{\mu}=\mu(|y|+|z|)$, then it is a $g$--expectation. This
domination inequality (\ref{domination}) has been applied as a testing
criteria in our data analysis. The results strongly support that both the
market pricing mechanism and the bid--ask pricing mechanism under our
investigation can be modelled as $g$--expectations, and that the bid--ask
prices are then produced by this single mechanism.

The paper is organized as follows: in Section \ref{ss2} we present the notion
of $g$--expectation and show that, for each well--defined function $g$ it
satisfies the basic conditions (A1)--(A4) of a dynamic pricing mechanism of
derivatives. We then show that, a dynamic pricing mechanism dominated by a
$g_{\mu}$--expectation, i.e., (\ref{domination}) is satisfied, is a
$g$--expectation. In Section 3, we will present some equivalent conditions to
show that the behaviors of a $g$--expectation are perfectly reflected by its
generating function $g$. We also provide some examples and explain how to
statistically find the function $g$ by testing the input--output data of prices.

In Appendix 4.2 we apply the crucial domination inequality (\ref{domination})
to test the market pricing mechanisms and the bid--ask pricing mechanisms of
S\&P500 index future options and S\&P500 index options, using data of
parameter files provided by CME and CBOE. The result supports that they are
$g$--expectations.

Application of the dynamic expectations and pricing mechanisms is to risk
measures. Axiomatic conditions for a (one step) coherent risk measure was
introduced by Artzner, Delbaen, Eber and Heath 1999 \cite{ADEH1999} and, for a
convex risk measure, by F\"{o}llmer and Schied (2002) {\cite{Fo-Sc}}. Rosazza
Gianin (2003) studied dynamic risk measures using the notion of $g$%
--expectations in \cite{rosazza} (see also \cite{Peng2003b}, \cite{El-Bar},
\cite{El-Bar2005}) in which an additional condition of cash translatability is assumed.

\section{The pricing mechanisms and $g$--pricing mechanisms by BSDE\label{ss2}%
}

Let us consider a market of financial derivatives in which the price
$(S_{t})_{t\geq0}$ of the underlying assets is driven by a $d$--dimensional
Brownian motion $(B_{t})_{t\geq0}$ in a probability space $(\Omega
,\mathcal{F},P)$. Here $S$ is an $m$--dimensional process, namely the number
of the underlying assets is $m$. We assume that the past information
$\mathcal{F}_{t}^{S}$ of the price $S$ before $t$ coincides with that of the
Brownian motion:
\[
\mathcal{F}_{t}^{S}=\sigma \{S_{s},\;s\leq t\}=\mathcal{F}_{t}:=\sigma
\{B_{s},\;s\leq t\}.
\]
A derivative $X$ with maturity $T$ is an $\mathcal{F}_{T}$--measurable and
square--integrable random value called maturity value is denoted by $X\in
L^{2}(\mathcal{F}_{T})$. The market price $Y_{t}$ of this derivative at time
$t<T$ is assumed to be in $L^{2}(\mathcal{F}_{t})$.

Let us consider a BSDE model of a pricing mechanism of derivatives, where
$Y_{t}$ is the solution of the following BSDE:%

\begin{equation}
Y_{t}=X+\int_{t}^{T}g(s,Y_{s},Z_{s})ds-\int_{t}^{T}Z_{s}dB_{s}.\label{tsBSDE}%
\end{equation}
Here $(Y,Z)$ a pair of the adapted processes to be solved, $g$ is a given
function
\[
g:(\omega,t,y,z)\in \Omega \times \lbrack0,\infty)\times R\times R^{d}\rightarrow
R.
\]
We call $g$ the generating function of the BSDE. It satisfies the following
basic assumptions for each $\forall y$, $\bar{y}\in R$ and $z$, $\bar{z}\in
R^{d}$,
\begin{equation}
\left \{
\begin{array}
[c]{l}%
g(\cdot,y,z)\in L_{\mathcal{F}}^{2}(0,T),\  \  \forall T\in(0,\infty),\\
|g(t,y,z)-g(T,\bar{y},\bar{z})|\leq \mu(|y-\bar{y}|+|z-\bar{z}|).
\end{array}
\right.  \label{h2.1}%
\end{equation}
It is important to consider the following special situation:
\begin{equation}
\left \{
\begin{array}
[c]{rrl}%
\text{(a)\  \  \ } & g(\cdot,0,0) & \equiv0,\\
\text{(b)\  \  \ } & g(\cdot,y,0) & \equiv0,\; \forall y\in R.
\end{array}
\  \right.  \label{h2.2}%
\end{equation}
Obviously (b) implies (a). This BSDE (\ref{tsBSDE}) was introduced by Bismut
\cite{Bismut73}, \cite{Bismut78} for the case where $g$ is a linear function
of $(y,z)$. \cite[Pardoux-Peng, 1990]{Pardoux-Peng1990} obtained the following
fundamental result: for each $X\in L^{2}(\mathcal{F}_{T})$, there exists a
unique square--integrable adapted solution $(Y,Z)$ of the BSDE (\ref{tsBSDE}).
The following notion of $g$--expectations was introduced by \cite[Peng
1997a]{Peng1997a} and \cite[Peng 1997]{Peng1997}.

\begin{definition}
\label{EgstX}We denote by $\mathbb{E}_{t,T}^{g}[X]:=Y_{t}$:%
\begin{equation}
\mathbb{E}_{t,T}^{g}[\cdot]:L^{2}(\mathcal{F}_{T})\rightarrow L^{2}%
(\mathcal{F}_{t}),\  \  \ 0\leq t\leq T<\infty. \label{Def2.2aa}%
\end{equation}
$(\mathbb{E}_{t,T}^{g}[\cdot])_{0\leq t\leq T<\infty}$ is called
$g$--expectation, or $g$--pricing mechanism.
\end{definition}

As an example, we consider the following Black--Scholes pricing mechanism:

\begin{example}
(\textbf{Black--Scholes Pricing Mechanism}) Consider a financial market
consisting of $2$ underlying assets: one bond and one stock. We denote by
$S_{0}(t)$ the price of the bond and by $S(t)$ the price of the stock at time
$t$. We assume that $S_{0}(t)$ satisfies an ordinary differential equation:
$dS_{0}(t)=r_{t}S_{0}(t)dt,$ and $S(t)$ is the solution of the following
stochastic differential equation (SDE) with $1$--dimensional Brownian motion
$B$ (i.e., $d=1$) as driven noise:
\[
dS(t)=S(t)(b_{t}dt+\sigma_{t}dB_{t}),\  \ S(0)=p.
\]
Here $r_{t}$ is the interest rate, $b_{t}$ the rate of the expected return and
$\sigma_{t}$ the volatility of the stock at the time $t$. $r_{t}$, $b_{t}$,
$\sigma_{t}$ and $\sigma^{-1}$ are assumed to be $\mathcal{F}_{t}$--measurable
and uniformly bounded. Black and Scholes have solved the problem of the market
pricing mechanism of an European call option $X=(S_{T}-k)^{+}$ and put option
$X=(k-S_{T})^{+}$, where $k$ is the strike price, under the assumption that
$r$, $b$ and $\sigma$ are constant. Their main idea can be easily adapted to
our slightly more general situation for a derivative $X\in L^{2}%
(\mathcal{F}_{T})$ with maturity $T$. Consider an investor with the following
investment portfolio at a time $t\leq T$: he has $n_{0}(t)$ bonds and $n(t)$
stock, i.e., he invests $n_{0}(t)S_{0}(t)$ in bond and $\pi(t)=n(t)S(t)$ in
the stock. We define by $Y_{t}$ the investor's wealth invested in the market
at time $t$:
\[
Y_{t}=n_{0}(t)P_{0}(t)+n(t)P(t).
\]
We make the so called \textquotedblleft self--financing
assumption\textquotedblright:
\[
dY_{t}=n_{0}(t)dS_{0}(t)+n(t)dS(t)
\]
or
\[
dY_{t}=[r_{t}Y_{t}+(b_{t}-r_{t})\pi(t)]dt+\sigma_{t}\pi_{t}dB_{t}.
\]
We denote $g(t,y,z):=-r_{t}y-(b_{t}-r_{t})\sigma_{t}^{-1}z$. Then, by denoting
$Z_{t}=\sigma_{t}\pi(t)$, the above equation is
\[
-dY_{t}=g(t,Y_{t},Z_{t})dt-Z_{t}dB_{t}.
\]
We observe that the above function $g$ satisfies (\ref{h2.1}). It follows from
the existence and uniqueness theorem of BSDE that for each derivative $X\in
L^{2}(\mathcal{F}_{T})$, there exists a unique adapted solution $(Y,Z)$ with
the terminal condition $Y_{T}=X$. This result of existence and uniqueness is
economically meaningful: in order to replicate the derivative $X$ at the
maturity $T$, the investor needs and only needs to invest the $Y_{t}$ at the
present time $t$ and then, during the time interval $s\in \lbrack t,T]$, to
perform the portfolio strategy $\pi(s)=\sigma_{s}^{-1}Z_{s}$. Furthermore, by
Comparison Theorem of BSDE, if he wants to replicate a derivative $\bar{X}$
with the same maturity $T$ which is bigger than $X$ (i.e., $\bar{X}\geq X$ and
$P(\bar{X}\geq X)>0$) then he must invest more than $Y_{t}$ $\ $at the time
$t$. This means there this no arbitrage opportunity. In this situation
$Y_{t}=\mathbb{E}_{t,T}^{g}[X]$ is called the Black--Scholes price, and
$(\mathbb{E}_{t,T}^{g}[\cdot])_{0\leq t\leq T<\infty}$ is called the
corresponding Black--Scholes pricing mechanism. We observe that the generating
function $g$ satisfies (a) of condition (\ref{h2.2}).
\end{example}

\begin{example}
The following problem was considered in \cite{Bergman}, \cite{CK1993} and
\cite{EPQ1997}: the investor is allowed to borrow money at time $t$ at an
interest rate $R_{t}>r_{t}$. The amount borrowed at time $t$ is equal to
$(Y_{t}-\pi(t))^{-}$. In this case the wealth process $Y$ still satisfies
BSDE:
\[
-dY_{t}=g(t,Y_{t},Z_{t})dt\,-\,Z_{t}dW_{t}.
\]
with $g(t,y,z):=-r_{t}y-(b_{t}-r_{t})\sigma_{t}^{-1}z+(R_{t}-r_{t}%
)(y-\sigma_{t}^{-1}z)^{-}$. This derives a $g$--pricing mechanism with a
sub--additive generating function $g$.
\end{example}

Similar equations appear in continuous trading with short sales constraints
with different risk premium for long and short positions (cf. \cite{JK1995},
\cite{HP1991} and \cite{EPQ1997}). In this case $g(t,y,z):=-r_{t}%
y-(b_{t}-r_{t})\sigma_{t}^{-1}z+k_{t}z^{-}$. We observe that in each of the
above three examples, $g$ is sub-additive in $(y,z)$.

The following result, obtained in \cite{Peng2003b}--Theorem 3.4, explains why
this $g$--expectation is a good candidate to model a dynamic pricing mechanism
of derivatives:

\begin{proposition}
\label{p2.1aa}Let the generating function $g$ satisfies (\ref{h2.1}) and
(\ref{h2.2})--(a). Then the above defined $g$--expectation $\mathbb{E}%
^{g}[\cdot]$ is a dynamic pricing mechanism of derivatives, i.e., it
satisfies, for each $t\leq T<\infty$, $X$, $\bar{X}\in L^{2}(\mathcal{F}_{T}%
)$, \newline \textbf{(A1)} $\mathbb{E}_{t,T}^{g}[X]\geq \mathbb{E}_{t,T}%
^{g}[\bar{X}]$, a.s., if $X\geq \bar{X}$;\newline \textbf{(A2)} $\mathbb{E}%
_{T,T}^{g}[X]=X$; \newline \textbf{(A3)} $\mathbb{E}_{s,t}^{g}[\mathbb{E}%
_{t,T}^{g}[X]]=\mathbb{E}_{s,T}^{g}[X]$; for $s\leq t$;\newline \textbf{(A4)}
$1_{A}\mathbb{E}_{t,T}^{g}[X]=\mathbb{E}_{t,T}^{g}[1_{A}X]$,\ $\forall
A\in \mathcal{F}_{t}$,\newline where $I_{A}$ is the indicator function of $A$,
i.e., $I_{A}(\omega)$ equals to $1$, when $\omega \in A$ and $0$ otherwise.
\newline
\end{proposition}

\begin{remark}
(A1) and (A2) are economically obvious conditions for a pricing mechanism.
Condition (A3) means that, at the time $s$, the random value $\mathbb{E}%
_{t,T}^{g}[X]$ can be regarded as a maturity value with maturity $t$. The
price of this derivative at $s$ is $\mathbb{E}_{s,t}^{g}[\mathbb{E}_{t,T}%
^{g}[X]]$. It must be the same as the price $\mathbb{E}_{s,T}^{g}[X]$ of $X$
at $s$.
\end{remark}

\begin{remark}
The meaning of condition (A4) is that, since at time $t$, the agent knows the
value of whether $I_{A}$ is $1$ or $0$. \ When $I_{A}$ is $1$, then the price
$\mathbb{E}_{t,T}^{g}[1_{A}X]$ of $1_{A}X$ must be the same as $\mathbb{E}%
_{t,T}^{g}[X]$, otherwise $1_{A}X=0$, so it worthies $0$.
\end{remark}

From the above results we see that $\mathbb{E}^{g}$ is a good candidate to be
a dynamic pricing mechanism. The following result provides a criteria to test
if a dynamic pricing mechanism is a $g$--expectation. The proof can be found
in \cite{Peng2005d}.

\begin{definition}
A system of mappings $(\mathbb{E}_{t,T}[\cdot])_{0\leq t\leq T<\infty}$
\[
\mathbb{E}_{t,T}[X]:X\in L^{2}(\mathcal{F}_{T})\mapsto L^{2}(\mathcal{F}_{t})
\]
is called a dynamic pricing mechanism of derivatives if it satisfies
(A1)--(A4) (with $\mathbb{E}\mathbf{[\cdot]}$ in the place of $\mathbb{E}%
^{g}[\cdot]$).
\end{definition}

\begin{theorem}
\label{m7.1}Let $\mathbb{E}_{t,T}[\cdot]_{0\leq t\leq T<\infty}$ be an dynamic
pricing mechanism. If there exists a sufficiently large constant $\mu>0$, such
that the following domination criteria is satisfied\newline \textbf{(A5):}\
\begin{equation}
\mathbb{E}_{t,T}[X]-\mathbb{E}_{t,T}[\bar{X}]\leq \mathbb{E}_{t,T}^{g_{\mu}%
}[X-\bar{X}]. \label{domination}%
\end{equation}
$\mathbb{E}^{g_{\mu}}$ is a $g$--expectation with the generating function
$g_{\mu}$ defined by
\begin{equation}
g_{\mu}(y,z):=\mu|y|+\mu|z|,\;(y,z)\in R\times R^{d}. \label{e2.9}%
\end{equation}
Then there exists a unique generating function $g(\omega,t,y,z)$ satisfying
(\ref{h2.1}) and (\ref{h2.2})-(a) such that, for each $t\leq T$ and for each
derivative $X\in L^{2}(\mathcal{F}_{T})$, we have
\begin{equation}
\mathbb{E}_{t,T}[X]=\mathbb{E}_{t,T}^{g}[X], \label{e7.1}%
\end{equation}
namely $\mathbb{E}$ is a $g$--expectation.
\end{theorem}

\begin{remark}
This theorem also implies that, for a generating function $g$ satisfying
(\ref{h2.1}) and (\ref{h2.2})--(a), the corresponding $g$--expectation
$\mathbb{E}^{g}$ is also dominated by $\mathbb{E}^{g_{\mu}}$, i.e., (A5) is
satisfied. This can be also directly proved by using the comparison theorem of BSDE.
\end{remark}

\begin{remark}
It turns out that the domination condition (\ref{domination}) becomes a
crucial criteria to test whether a dynamic pricing mechanism of derivatives is
a $g$--expectation. We provide a test in Appendix 4.2 to use market data to
check the inequality (\ref{domination}).
\end{remark}

\begin{remark}
This deep result has non-trivially generalized the main result of
\cite{CHMP2002}, theoretically and practically, where a special case
$g=g(t,z)$ with $g(s,0)\equiv0$ is considered. The $g$--expectation originally
introduced in \cite{Peng1997} corresponds such situation of \textquotedblleft
zero interest rate\textquotedblright. (cf. Proposition \ref{m2a4}, or
\cite{Peng2003b}).
\end{remark}

\textbf{Markovian pricing mechanisms }We limit ourselves to consider, for each
fixed maturity $T$, the derivatives $X$ depending only on the price $S_{T}$,
i.e., $X$ is a path independent derivative. $X$ is then in the class of
\[
X=\Phi(S_{T})\text{ with }\Phi \in L^{2}(S_{T})
\]
where $L^{2}(S_{T})$ denotes the collection of all real functions $\Phi$
defined on $R^{n}$ such that $\Phi(S_{T})\in L^{2}(\mathcal{F}_{T})$. A
dynamic pricing mechanism $\mathbb{E}$ is called a Markovian pricing mechanism
if for each $0\leq t\leq T<\infty$ and $\Phi \in L^{2}(S_{T})$ there exists
$\phi \in L^{2}(S_{t})$ such that $\mathbb{E}_{t,T}[\Phi(S_{T})]=\phi(S_{t})$.
In other words, \emph{the price of a path--independent option by a Markovian
pricing mechanism is still path-independent. }

\begin{example}
We consider a situation where the underlying price $S$ is a diffusion process:%
\[
dS_{t}=b(S_{t})dt+\Lambda(S_{t})dB_{t},\  \ S_{0}=s_{0}\in R^{n}.
\]
where $b$ and $\Lambda$ are given Lipschitz functions of $R^{n}$ valued on
$R^{n}$ and $R^{n\times d}$ respectively. If a generating function $g$ has the
following form:
\[
g(t,y,z)=f(S_{t},y,z),
\]
where $f$ is a Lipschitz function of $(s,y,z)\in R^{n}\times R\times R^{d}$.
By the nonlinear Feynman--Kac formula introduced in \cite[Peng 1991]%
{Peng1991}, \cite[Peng 1992]{Peng1992} and developed in \cite[Pardoux-Peng
1992]{PP1992}, for each option $X=\Phi(S_{T})$ with smooth function $\Phi$ the
price of the related $g$--expectation is
\[
\mathbb{E}_{t,T}^{g}[\Phi(S_{T})]=u(t,S_{t})
\]
where $u:R_{+}\times R^{n}\longmapsto R$ is the (viscosity) solution of the
following PDE defined on $(t,s)\in \lbrack0,T]\times R^{n}$:%
\[
\frac{\partial u}{\partial t}+\frac{1}{2}\sum_{i,j=1}^{n}(\Lambda
(s)\Lambda^{T}(s))_{ij}\frac{\partial^{2}u}{\partial s_{i}\partial s_{j}}%
+\sum_{i=1}^{n}b_{i}(s)\frac{\partial u}{\partial s_{i}}+f(s,u,\Lambda
^{T}(s)\nabla u)=0
\]
with terminal condition $u(T,s)=\Phi(s)$. If $S_{t}$ is a $1$--dimensional
geometric Brownian motion, i.e., $\Lambda(s)=\sigma s$ and $b(s)=\mu s$, then
the above PDE becomes%
\[
\frac{\partial u}{\partial t}+\frac{1}{2}\sigma^{2}s^{2}\frac{\partial^{2}%
u}{\partial s^{2}}+\mu s\frac{\partial u}{\partial s}+f(s,u,\sigma
s\frac{\partial u}{\partial s})=0.
\]
The Black--Scholes formula corresponds to the case $f=-ry-(\mu-r)\sigma^{-1}%
z$. We then have%
\[
\frac{\partial u}{\partial t}+\frac{1}{2}\sigma^{2}s^{2}\frac{\partial^{2}%
u}{\partial s^{2}}-ru+rs\frac{\partial u}{\partial s}=0,\  \ u(T,s)=\Phi(s).
\]

\end{example}

\section{Characterization of $g$-pricing mechanism by its generating function
$g$}

For a pricing mechanism, it is important to distinguish the selling price and
buying price of a same pricing mechanism, corresponding to the ask price and
the bid price if the mechanism under investigation is generated through the
system of market makers of an option market (cf. \cite{Hull} Sec. 6.5 and
\cite{Musiela-Rutkowski}). If $\mathbb{E}_{t,T}[X]$ is the ask price at the
time $t$ of a derivative $X$ with maturity $T$, then the bid price must be
$-\mathbb{E}_{t,T}[-X]$ and we have, in general,
\[
\mathbb{E}_{t,T}[X]>-\mathbb{E}_{t,T}[-X].
\]
Here we stress our point of view that, in fact, the ask price and bid price
are produced by a single mechanism, called bid--ask pricing mechanism of
market makers. Our result of data analysis to test the criteria (A5) of the
domination condition (\ref{domination}) strongly supports this point of view.
Moreover, this analysis also supports our point of view that, for a
well--developed market, there exist a function $g$ satisfying Lipschitz
condition (\ref{h2.1}) such that the corresponding ask-bid pricing mechanism
is modeled by the $g$--expectation $\mathbb{E}^{g}[\cdot]$. 

A rational dynamic pricing mechanism also possesses some other important
properties, such as convexity, sub-additivity. See \cite{ADEH1999},
\cite{El-Bar}, \cite{El-Bar2005}, \cite{Chen-Epstein2002},
\cite{ElKaroui-Quenez}, \cite{EK-Q}, \cite{EPQ1997}, \cite{Fo-Sc},
\cite{Fritteli}, \cite{Peng-Yang}, \cite{rosazza}, \cite{JiangL2004a},
\cite{JiangL2004b} among many others. We will see that the generating
function$\ g$ perfectly reflects the behavior of $\mathbb{E}^{g}$. This may
play an important role to statistically find $g$ by using the corresponding
data of prices. \ In the following we provide several theoretical results with
proofs given in Appendix. This problem was treated also by \cite{rosazza},
\cite{JiangL2004a} and \cite{JiangL2004b}.

\begin{proposition}
\label{invcompar}Let $g$, $\bar{g}:$ $(\omega,t,y,z)\in \Omega \times
\lbrack0,\infty)\times R\times R^{d}\rightarrow R$ be two generating functions
satisfying (\ref{h2.1}). Then the following two conditions are equivalent:
\newline \textbf{(i) }$g(\omega,t,y,z)\geq \bar{g}(\omega,t,y,z)$,
$\forall(y,z)\in R\times R^{d}$, $dP\times dt\ $a.s. \newline \textbf{(ii) }The
corresponding $g$--pricing mechanisms $\mathbb{E}^{g}[\cdot]$ and
$\mathbb{E}^{\bar{g}}[\cdot]$ satisfy
\[
\mathbb{E}_{t,T}^{g}[X]\geq \mathbb{E}_{t,T}^{\bar{g}}[X],\  \forall0\leq t\leq
T<\infty,\  \forall X\in L^{2}(\mathcal{F}_{T}).\
\]
In particular, $\mathbb{E}^{g}[X]\equiv \mathbb{E}^{\bar{g}}[X]$ if and only if
$g\equiv \bar{g}$.
\end{proposition}

\begin{corollary}
\label{Maker-E-X}The following two conditions are equivalent: \newline%
\textbf{(i)} The generating function $g$ satisfies, for each $(y,z)\in R\times
R^{d}$,
\[
g(t,y,z)\geq-g(t,-y,-z),\; \text{a.e., a.s.,}%
\]
\textbf{(ii)} $\mathbb{E}_{t,T}^{g}[\cdot]:L^{2}(\mathcal{F}_{T})\longmapsto
L^{2}(\mathcal{F}_{t})$ satisfies, for each $0\leq t\leq T$,
\[
\mathbb{E}_{t,T}^{g}[X]\geq-\mathbb{E}_{t,T}^{g}[-X],\  \ X\in L^{2}%
(\mathcal{F}_{T}).
\]

\end{corollary}

\begin{proof}
We denote $\bar{g}(t,y,z):=-g(t,-y,-z)$ and compare the following two BSDE:%
\[
Y_{t}=X+\int_{t}^{T}g(s,Y_{s},Z_{s})ds-\int_{t}^{T}Z_{s}dB_{s},\ t\in
\lbrack0,T],
\]
and
\[
\bar{Y}_{t}=-X+\int_{t}^{T}g(s,\bar{Y}_{s},\bar{Z}_{s})ds-\int_{t}^{T}\bar
{Z}_{s}dB_{s},\ t\in \lbrack0,T].
\]
or, with $\hat{Y}=-\bar{Y}$, $\hat{Z}=-\bar{Z}$%
\[
\hat{Y}_{t}=X+\int_{t}^{T}\bar{g}(s,\hat{Y}_{s},\hat{Z}_{s})ds-\int_{t}%
^{T}\hat{Z}_{s}dB_{s},\ t\in \lbrack0,T].
\]
From the above Proposition it follows that $\mathbb{E}_{t,T}^{g}[\cdot
]\geq \mathbb{E}_{t,T}^{\bar{g}}[\cdot]$ iff $g\geq \bar{g}$. This with
$\mathbb{E}_{t,T}^{\bar{g}}[X]=-\mathbb{E}_{t,T}^{g}[-X]$ yields (i)
$\Leftrightarrow$ (ii).
\end{proof}

\begin{proposition}
\label{Proposition-2}The following two conditions are equivalent:
\newline \textbf{(i)} The generating function $g=g(t,y,z)$ is convex (resp.
concave) in $(y,z)$, i.e., for each $(y,z)$ and $(\bar{y},\bar{z})$ in
$R\times R^{d}$ and for a.e. $t\in \lbrack0,T]$,
\begin{align*}
g(t,\alpha y+(1-\alpha)\bar{y},\alpha z+(1-\alpha)\bar{z})  &  \leq \alpha
g(t,y,z)+(1-\alpha)g(t,\bar{y},\bar{z}),\; \text{a.s.}\\
\; \; \text{(resp.}  &  \geq \alpha g(t,y,z)+(1-\alpha)g(t,\bar{y},\bar{z}),\;
\text{a.s.).}%
\end{align*}
\textbf{(ii)} The corresponding pricing mechanism $(\mathbb{E}_{t,T}^{g}%
[\cdot])_{0\leq t\leq T<\infty}$ is convex (resp. concave), i.e., for each
fixed $\alpha \in \lbrack0,1]$, we have
\begin{align}
\mathbb{E}_{t,T}^{g}[\alpha X+(1-\alpha)\bar{X}]  &  \leq \alpha \mathbb{E}%
_{t,T}^{g}[X]+(1-\alpha)\mathbb{E}_{t,T}^{g}[\bar{X}],\; \text{a.s.}%
\label{Econvex}\\
\text{(resp. }  &  \geq \alpha \mathbb{E}_{t,T}^{g}[X]+(1-\alpha)\mathbb{E}%
_{t,T}^{g}[\bar{X}],\; \text{a.s.)}\nonumber \\
\text{for each }t  &  \leq T,\  \text{and }X,\bar{X}\in L^{2}(\mathcal{F}%
_{T}).\nonumber
\end{align}

\end{proposition}

\begin{proposition}
\label{Proposition-3}The following two conditions are equivalent:
\newline \textbf{(i)} The generating function $g$ is positively homogenous in
$(y,z)\in R\times R^{d}$, i.e.,
\[
g(t,\lambda y,\lambda z)=\lambda g(t,y,z),\; \text{a.e., a.s.,}%
\]
\textbf{(ii)} The corresponding pricing mechanism $\mathbb{E}_{t,T}^{g}%
[\cdot]:L^{2}(\mathcal{F}_{T})\longmapsto L^{2}(\mathcal{F}_{t})$ is
positively homogenous: for each $0\leq t\leq T$ , i.e., $\mathbb{E}_{t,T}%
^{g}[\lambda X]=\lambda \mathbb{E}_{t,T}^{g}[X]$, for each $\lambda \geq0$ and
$X\in L^{2}(\mathcal{F}_{T})$.
\end{proposition}

From the above two propositions we immediately have

\begin{corollary}
The following two conditions are equivalent: \newline \textbf{(i)} The
generating function $g$ is sub-additive: for each $(y,z),\ (\bar{y},\bar
{z})\in R\times R^{d}$,
\[
g(\omega,t,y+\bar{y},z+\bar{z})\leq g(\omega,t,y,z)+g(\omega,t,\bar{y},\bar
{z}),\;dt\times dP\text{, a.s.,}%
\]
\textbf{(ii)} The corresponding pricing mechanism $\mathbb{E}_{t,T}^{g}%
[\cdot]:L^{2}(\mathcal{F}_{T})\longmapsto L^{2}(\mathcal{F}_{t})$ is is
sub-additive: for each $0\leq t\leq T$ and $X$, $\bar{X}\in L^{2}%
(\mathcal{F}_{T})$
\[
\mathbb{E}_{t,T}^{g}[X+\bar{X}]\leq \mathbb{E}_{t,T}^{g}[X]+\mathbb{E}%
_{t,T}^{g}[\bar{X}].
\]

\end{corollary}

\begin{proposition}
\label{Proposition-4}The generating function $g$ is independent of $y$ if and
only if the corresponding $g$--expectation satisfies the following
\textquotedblleft cash translatability\textquotedblright \ property: for each
$t\leq T$,
\[
\mathbb{E}_{t,T}^{g}[X+\eta]=\mathbb{E}_{t,T}^{g}[X]+\eta,\; \text{a.s.},\text{
for each }X\in L^{2}(\mathcal{F}_{T}),\  \eta \in L^{2}(\mathcal{F}_{t}).
\]

\end{proposition}

We consider the following self--financing condition:
\[
\mathbb{E}_{t,T}^{g}[0]\equiv0,\; \forall0\leq t\leq T.
\]

\begin{proposition}
$\mathbb{E}^{g}[\cdot]$ satisfies the self--financing condition if and only if
its generating function $g$ satisfies (\ref{h2.2})--(a).
\end{proposition}

\begin{proof}
The \textquotedblleft if\textquotedblright \ part is obvious. \newline The
\textquotedblleft only if part\textquotedblright: $Y_{t}:=\mathbb{E}_{t,T}%
^{g}[0]\equiv0$, implies
\[
Y_{t}\equiv0\equiv0+\int_{t}^{T}g(s,0,Z_{s})ds-\int_{t}^{T}Z_{s}dB_{s}%
,\;t\in \lbrack0,T].
\]
Thus $Z_{t}\equiv0$ and then $g(t,0,Z_{t})=g(t,0,0)\equiv0$.
\end{proof}

\medskip

\textquotedblleft Zero--interest rate\textquotedblright \ condition:
\[
\mathbb{E}_{t,T}^{g}[\eta]=\eta,\; \forall0\leq t\leq T<\infty \text{, }\eta \in
L^{2}(\mathcal{F}_{t}).
\]

\begin{proposition}
\label{m2a4}$\mathbb{E}^{g}[\cdot]$ satisfies the zero--interest rate
condition if and only if its generating function $g$ satisfies (\ref{h2.2})--(b).
\end{proposition}

\begin{proof}
\textbf{ }For a fixed $y\in R$, we consider $Y_{t}:=\mathbb{E}_{t,T}%
^{g}[y]\equiv y$. Let $Z_{t}$ be the corresponding It\^{o}'s integrand in
$Y$:
\[
Y_{t}=y+\int_{t}^{T}g(s,Y_{s},Z_{s})-\int_{t}^{T}Z_{s}dB_{s}\equiv y.
\]
But this is equivalent to
\[
Y_{t}\equiv y,\;Z_{s}\equiv0\text{, }g(s,y,0)\equiv0.
\]

\end{proof}

\medskip

For each $\bar{z}_{\cdot}^{i_{0}}\in L_{\mathcal{F}}^{2}(0,T)$%
\begin{equation}
\mathbb{E}_{t,T}[X]+\int_{0}^{t}\bar{z}_{s}^{i_{0}}dB_{s}^{i_{0}}%
=\mathbb{E}_{t,T}[X+\int_{t}^{T}\bar{z}_{s}^{i_{0}}dB_{s}^{i_{0}}]
\label{zBi0}%
\end{equation}

\begin{proposition}
\label{Proposition-5}Condition (\ref{zBi0}) holds if and only if $g(t,y,z)$
does not depends on the $i_{0}$--th component $z^{i_{0}}$ of $z\in R^{d}$.
\end{proposition}

\begin{proposition}
The following condition are equivalent: \newline \textbf{(i)} For each $0\leq
t\leq T$ and $X\in L^{2}(\mathcal{F}_{T}^{t})$, the $g$--pricing mechanism
$\mathbb{E}_{t,T}^{g}[X]$ is a deterministic number; \newline \textbf{(ii) }The
corresponding pricing generating function $g$ is a deterministic function of
$(t,y,z)\in \lbrack0,T]\times R\times R^{d}$. $\;$
\end{proposition}

The proof is similar as the others. We omit it.

\begin{example}
An interesting problem is: if we know that a pricing mechanism under our
investigation is a $g$--expectation $\mathbb{E}^{g}$, how to find the
generating function $g$? If we limited ourselves to only take data of prices
quoted by markets, this is still an open problem. We now consider a case of
\textquotedblleft toy model\textquotedblright \ where $g$ depends only on $z$,
i.e., $g=g(z):\mathbf{R}^{d}\rightarrow \mathbf{R}$. We will find such $g$ by
the following testing method. Let $\bar{z}\in R^{d}$ be given. We denote
$Y_{s}:=\mathbb{E}_{s,T}^{g}[\bar{z}\cdot(B_{T}-B_{t})]$, $s\in \lbrack t,T]$,
where $t$ is the present time. It is the solution of the following BSDE
\[
Y_{s}=\bar{z}\cdot(B_{T}-B_{t})+\int_{s}^{T}g(Z_{u})du-\int_{s}^{T}Z_{u}%
dB_{u},\;s\in \lbrack t,T].
\]
It is seen that the solution is $Y_{s}=\bar{z}\cdot(B_{s}-B_{t})+\int_{s}%
^{T}g(\bar{z})ds$, $Z_{s}\equiv \bar{z}$. Thus
\[
\mathbb{E}_{t,T}^{g}[\bar{z}\cdot(B_{T}-B_{t})]=Y_{t}=g(\bar{z})(T-t),
\]
or
\begin{equation}
g(\bar{z})=(T-t)^{-1}\mathbb{E}_{t,T}^{g}[\bar{z}\cdot(B_{T}-B_{t}%
)].\label{e2.exm1}%
\end{equation}
Thus the function $g$ can be tested as follows: at the present time $t$: if
the valuation $\mathbb{E}_{t,T}^{g}[\bar{z}\cdot(B_{T}-B_{t})]$ of (a toy
model of) derivative $\bar{z}\cdot(B_{T}-B_{t})$ is obtained, then $g(\bar
{z})$ is explicitly given by (\ref{e2.exm1}). We observe that, in the case
where $S$ is a geometric Brownian motion, $B_{T}-B_{t}$ can be expressed as a
function of $S_{T}/S_{t}$. But this cannot be applied to a general situation.
\end{example}

\begin{remark}
The above test is also applied for the case $g:[0,\infty)\times \mathbf{R}%
^{d}\rightarrow \mathbf{R}$, or for a more general situation $g=\gamma
y+g_{0}(t,z)$.
\end{remark}

An interesting problem is, in general, how to find the generating function $g$
by a testing of the input--output behavior of $\mathbb{E}^{g}[\cdot]$? Let
$b:R^{n}\longmapsto R^{n}$, $\bar{\sigma}:R^{n}\longmapsto R^{n\times d}$ be
two Lipschitz functions. For each $(t,x)\in R_{+}\times R^{n}$, we consider
the SDE of the form
\[
X_{s}^{t,x}=x+\int_{t}^{s}b(X_{s}^{t,x})ds+\int_{t}^{s}\sigma(X_{s}%
^{t,x})dB_{s},\; \;s\geq t.
\]
This SDE is regarded as the equation of the price of the underlying stock. The
following result was obtained in Proposition 2.3 of \cite{BCHMP00}.

\begin{proposition}
\label{Prop-Rep}We assume that the generating function $g$ satisfies
(\ref{h2.1})\textbf{. }We also assume that, for each fixed $(y,z)$,
$g(\cdot,y,z)\in D_{\mathcal{F}}^{2}(0,T)$ (the space of all $\mathcal{F}_{t}%
$--adapted processes with RCLL paths). Then for each $(t,x,p,y)\in
\lbrack0,\infty)\times R^{n}\times R^{n}\times R$, we have
\[
L^{2}\hbox{--}\lim_{\epsilon \rightarrow0}\frac{1}{\epsilon}[\mathbb{E}%
_{t,t+\epsilon}^{g}[y+p\cdot(X_{t+\epsilon}^{t,x}-x)]-y]=g(t,y,\sigma
^{T}(x)p)+p\cdot b(x).\label{limE}%
\]

\end{proposition}

\section{Appendix \label{ss8}}

\subsection{Proofs of Propositions 3.1--3.10}

We begin with introducing some technique lemmas. The first one is called
decomposition theorem of $\mathbb{E}^{g}$--supermartingale. The proof can be
find in \cite{Peng1999} and \cite{Peng2003b}.

\begin{proposition}
\label{p2.3} We assume (\ref{h2.1}). Let $Y\in D_{\mathcal{F}}^{2}(0,T)$ be an
$\mathbb{E}^{g}$--supermartingale, namely, for each $0\leq s\leq t\leq T$,%
\[
\mathbb{E}_{s,t}^{g}[Y_{t}]\leq Y_{t}\text{.}%
\]
Then there exists a unique $\mathcal{F}_{t}$--adapted increasing and RCLL
process $A\in D_{\mathcal{F}}^{2}(0,T)$ (thus predictable) with $A_{0}=0$,
such that, $Y$ is the solution of the following BSDE:%
\[
Y_{t}=Y_{T}+(A_{T}-A_{t})+\int_{t}^{T}g(s,Y_{s},Z_{s})ds-\int_{t}^{T}%
Z_{s}dB_{s},\  \ t\in \lbrack0,T].
\]

\end{proposition}

Let a function $f:(\omega,t,y,z)\in \Omega \times \lbrack0,T]\times R\times
R^{d}\rightarrow R$ satisfy the same Lipschitz condition (\ref{h2.1}) as for
$g$. For each fixed $(t,y,z)\in \lbrack0,T]\times R\times R^{d}$, we consider
the following SDE of It\^{o}'s type defined on $[t,T]$:
\begin{equation}
Y_{s}^{t,y,z}=y-\int_{t}^{s}f(r,Y_{r}^{t,y,z},z)dr+z\cdot(B_{s}-B_{t})
\label{Ytyz}%
\end{equation}
We have the following classical result of It\^{o}'s SDE.

\begin{lemma}
\label{em7.2a}We assume that $f$ satisfies the same Lipschitz condition
(\ref{h2.1}) as for $g$. Then there exists a constant $C$, depending only on
$\mu$, $T$ and $E\int_{0}^{T}|f(\cdot,0,0)|^{2}ds$, such that, for each
$(t,y,z)\in \lbrack0,T]\times R\times R^{d}$, we have
\begin{equation}
E[|Y_{s}^{t,y,z}-y|^{2}]\leq C(|y|^{2}+|z|^{2}+1)(s-t),\; \forall s\in \lbrack
t,T]. \label{ee7.15}%
\end{equation}

\end{lemma}

\begin{proof}
\textbf{ }It is classic that $E\int_{0}^{T}|f(r,Y_{r}^{t,y,z},z)|dr^{2}\leq
C_{0}(|y|^{2}+|z|^{2}+1)$, where $C_{0}$ depends only on $\mu$, $T$ and
$E\int_{0}^{T}|f(\cdot,0,0)|^{2}ds$. We then have
\begin{align*}
E[|Y_{s}^{t,y,z}-y|^{2}]  &  \leq2E[|\int_{t}^{s}f(r,Y_{r}^{t,y,z}%
,z)dr|^{2}]+2|z|^{2}(s-t)\\
\  &  \leq2E[\int_{t}^{s}|f(r,Y_{r}^{t,y,z},z)|^{2}dr]^{1/2}(t-s)+2|z|^{2}%
(s-t)\\
\  &  \leq C(|y|^{2}+|z|^{2}+1)(s-t).
\end{align*}

\end{proof}

For each $n=1,2,3,\cdots$, we set
\begin{align}
f^{n}(s,y,z)  &  :=\sum_{i=0}^{2^{n}-1}f(s,Y_{s}^{t_{i}^{n},y},z)1_{[t_{i}%
^{n},t_{i+1}^{n})}(s),\;s\in \lbrack0,T]\label{efnyz}\\
t_{i}^{n}  &  =i2^{-n}T,i=0,1,2,\cdots,2^{n}.
\end{align}
It is clear that $f^{n}$ is an $\mathcal{F}_{t}$--adapted process.

\begin{lemma}
\label{FntoF}For each fixed $(y,z)\in R\times R^{d}$, $\{f^{n}(\cdot
,y,z)\}_{n=1}^{\infty}$ converges to $f(\cdot,y,z)$ in $L_{\mathcal{F}}%
^{2}(0,T)$, i.e.,
\begin{equation}
\lim_{n\rightarrow \infty}E\int_{0}^{T}|f^{n}(s,y,z)-f(s,y,z)|^{2}ds=0.
\label{fntof}%
\end{equation}

\end{lemma}

\begin{proof}
For each $s\in \lbrack0,T)$, there are some integers $i\leq2^{n}-1$ such that
$s\in \lbrack t_{i}^{n},t_{i+1}^{n})$. We have, by (\ref{ee7.15})
\begin{align*}
E[|f^{n}(s,y,z)-f(s,y,z)|^{2}]  &  =E[|f(s,Y_{s}^{t_{i}^{n},y}%
,z)-f(s,y,z)|^{2}]\\
\  &  \leq \mu^{2}E[|Y_{s}^{t_{i}^{n},y,z}-y|^{2}]\\
\  &  \leq \mu^{2}C(|y|^{2}+|z|^{2}+1)2^{-n}T.
\end{align*}
Thus $\{f^{n}(\cdot,y,z)\}_{n=1}^{\infty}$ converges to $f(\cdot,y,z)$ in
$L_{\mathcal{F}}^{2}(0,T)$.
\end{proof}

\begin{lemma}
\label{Lem-f-f-}If for each $(t,y,z)\in \lbrack0,T]\times R\times R^{d}$, we
have
\[
f(\omega,r,Y_{r}^{t,y,z},z)\geq0\  \  \text{(resp. }=0\text{)},\;(\omega
,r)\in \lbrack t,T]\times \Omega \text{, .}dr\times dP\text{-a.s..}%
\]
Then, for each $(y,z)\in R\times R^{d}$,
\begin{equation}
f(\omega,t,y,z)\geq0,\; \text{(resp. }=0\text{)},\ (\omega,t)\in
\lbrack0,T]\times \Omega \text{, .}dt\times dP\text{-a.s..} \label{f-f-}%
\end{equation}

\end{lemma}

\begin{proof}
Let us fix $y$ and $z$. We define $f^{n}(s,y,z)$ as in (\ref{efnyz}). It is
clear that,
\[
f^{n}(r,y,z)\geq0,\  \text{(resp. }=0\text{),}\;(\omega,r)\in \lbrack
0,T]\times \Omega \text{, .}dr\times dP\text{ a.s.}%
\]
But from Lemma \ref{FntoF} we have $f^{n}(\cdot,y,z)\rightarrow f(\cdot,y,z)$,
in $L_{\mathcal{F}}^{2}(0,T)$ as $n\rightarrow \infty$. We thus have
\textbf{\ref{f-f-}}.
\end{proof}

\medskip

We now can give the proofs of several propositions given in the previous
section. The method is very different from \cite{rosazza}, \cite{JiangL2004a}
and \cite{JiangL2004b} where Proposition \ref{Prop-Rep} plays a central role.

\medskip

\textbf{Proof of Proposition \ref{invcompar} (i)}$\Rightarrow$(ii) is the
well--known comparison theorem of BSDE (cf. \cite{Peng1992} and \cite{EPQ1997}%
).\newline \textbf{ (ii)}$\Rightarrow$\textbf{(i): }For fixed $t\geq0$ and
$(y,z)$ in $R\times R^{d}$, let $Y^{t,y,z}$ be the solution of SDE
(\ref{Ytyz}) with $f=g$. From (ii) we have
\[
\mathbb{E}_{r,s}^{\bar{g}}[Y_{s}^{t,y,z}]\leq \mathbb{E}_{r,s}^{g}%
[Y_{s}^{t,y,z}]=Y_{s}^{t,y,z},\ t\leq r\leq s.
\]
Thus $(Y_{s}^{t,y,z})_{s\geq t}$ is an $\mathbb{E}^{\bar{g}}$%
--supermartingale. From the decomposition theorem, i.e., Proposition
\ref{p2.3}, it follows that there exists an increasing process $(\bar{A}%
_{s})_{s\geq t}$ such
\[
Y_{s}^{t,y,z}=y-\int_{t}^{s}\bar{g}(r,Y_{r}^{t,y,z},\bar{Z}_{r})dr-\bar{A}%
_{s}+\int_{t}^{s}\bar{Z}_{r}dB_{r},\ s\geq t.
\]
This with $Y_{s}^{t,y,z}=y-\int_{t}^{s}g(r,Y_{r}^{t,y,z},z)dr+\int_{t}%
^{s}zdB_{r}$ yields $\bar{Z}_{s}\equiv z$ and
\[
g(r,Y_{r}^{t,y,z},z)\geq \bar{g}(r,Y_{r}^{t,y,z},z),\ r\geq t.
\]
We then can apply the above Lemma \ref{Lem-f-f-} to prove that $g\geq \bar{g}$.
$\blacksquare$

\medskip

\textbf{Proof of Proposition \ref{Proposition-2} }We only prove the convex
case. \newline \textbf{(i)}$\Rightarrow$\textbf{(ii): }For a given $t>0$, we
set $Y_{s}^{X}:=\mathbb{E}_{s,t}^{g}[X]$, $Y_{s}^{\bar{X}}:=\mathbb{E}%
_{s,t}^{g}[\bar{X}]$, $s\in \lbrack0,t]$. These two pricing processes solve
respectively the following two BSDEs on $[0,t]$:
\begin{align*}
Y_{s}^{X} &  =X+\int_{s}^{t}g(r,Y_{r}^{X},Z_{r}^{X})dr-\int_{s}^{t}Z_{r}%
^{X}dB_{r},\\
Y_{s}^{\bar{X}} &  =\bar{X}+\int_{s}^{t}g(r,Y_{r}^{\bar{X}},Z_{r}^{\bar{X}%
})dr-\int_{s}^{t}Z_{r}^{\bar{X}}dB_{r}.
\end{align*}
Their convex combination: $(Y_{s},Z_{s}):=(\alpha Y_{s}^{X}+(1-\alpha
)Y_{s}^{\bar{X}},\alpha Z_{s}^{X}+(1-\alpha)Z_{s}^{\bar{X}})$, satisfies
\begin{align*}
Y_{s} &  =\alpha X+(1-\alpha)\bar{X}+\int_{s}^{t}[g(r,Y_{r},Z_{r})+\psi
_{r}]dr-\int_{s}^{t}Z_{r}dB_{r},\\
\text{where we set }\psi_{r} &  =\alpha g(r,Y_{r}^{X},Z_{r}^{X})+(1-\alpha
)g(r,Y_{r}^{\bar{X}},Z_{r}^{\bar{X}})-g(r,Y_{r},Z_{r})\text{.}%
\end{align*}
But since the price generating function $g$ is convex in $(y,z)$, we have
$\psi \geq0$. It then follows from the comparison theorem that $Y_{s}%
\geq \mathbb{E}_{s,t}^{g}[\alpha X+(1-\alpha)\bar{X}]$. We thus have
\textbf{(ii)}. \newline \textbf{(ii)}$\Rightarrow$\textbf{(i): }Let $Y^{t,y,z}$
be the solution of SDE (\ref{Ytyz}) with $f=g$. For fixed $t\in \lbrack0,T)$
and $(y,z)$, $(\bar{y},\bar{z})$ in $R\times R^{d}$, we have
\[
Y_{s}^{t,y,z}=\mathbb{E}_{r,s}^{g}[Y_{s}^{t,y,z}],\;Y_{s}^{t,\bar{y},\bar{z}%
}=\mathbb{E}_{r,s}^{g}[Y_{s}^{t,\bar{y},\bar{z}}],\  \ t\leq r\leq s.
\]
We set $Y_{s}:=\alpha Y_{s}^{t,y,z}+(1-\alpha)Y_{s}^{t,\bar{y},\bar{z}}$,
$s\in \lbrack t,T]$. By (\ref{Econvex}),
\begin{align*}
\mathbb{E}_{r,s}^{g}[Y_{s}] &  \leq \alpha \mathbb{E}_{r,s}^{g}[Y_{s}%
^{t,y,z}]+(1-\alpha)\mathbb{E}_{r,s}^{g}[Y_{s}^{t,\bar{y},\bar{z}}]\\
&  =\alpha Y_{r}^{t,y,z}+(1-\alpha)Y_{r}^{t,\bar{y},\bar{z}}=Y_{r}.
\end{align*}
Thus the process $Y$ is a $\mathbb{E}^{g}$--supermartingale defined on
$[t,T]$. It follows from the decomposition theorem, i.e., Proposition
\ref{p2.3}, that, there exists an increasing process $A$ such that
\[
Y_{s}=Y_{t}-\int_{t}^{s}g(r,Y_{r},Z_{r})dr-A_{s}+\int_{t}^{s}Z_{r}dB_{r}.
\]
We compare this with
\begin{align*}
Y_{s} &  =\alpha Y_{s}^{t,y,z}+(1-\alpha)Y_{s}^{t,\bar{y},\bar{z}}\\
&  =\alpha y+(1-\alpha)\bar{y}-\int_{t}^{s}[\alpha g(r,Y_{r}^{t,y,z}%
,z)+(1-\alpha)g(r,Y_{r}^{t,\bar{y},\bar{z}},\bar{z})]dr\\
&  +(\alpha z+(1-\alpha)\bar{z})\cdot(B_{s}-B_{t}),
\end{align*}
It follows that
\[
Y_{t}=\alpha y+(1-\alpha)\bar{y},\;Z_{r}\equiv \alpha z+(1-\alpha)\bar{z},\;
\]
Thus we have
\[
g(s,\alpha Y_{s}^{t,y,z}+(1-\alpha)Y_{s}^{t,\bar{y},\bar{z}},\alpha
z+(1-\alpha)\bar{z})\leq \alpha g(s,Y_{s}^{t,y,z},z)+(1-\alpha)g(s,Y_{s}%
^{t,\bar{y},\bar{z}},\bar{z}).
\]
We then can apply Lemma \ref{Lem-f-f-} to obtain (i). $\blacksquare$

\medskip

\textbf{Proof of Proposition \ref{Proposition-3}. }(i)$\Rightarrow$(ii) is
easy.\newline(ii)$\Rightarrow$(i):\textbf{ }Let $Y^{t,y,z}$ be the solution of
SDE (\ref{Ytyz}) with $f=g$. For fixed $t\in \lbrack0,T)$ and $(y,z)$ in
$R\times R^{d}$, we have $\lambda Y_{s}^{t,y,z}=\mathbb{E}_{s,t}^{g}[\lambda
Y_{t}^{t,y,z}]$, $s\in \lbrack t,T]$. This implies that, there exists a process
$Z_{\cdot}^{t,y,z,\lambda}$ such that
\[
\lambda Y_{s}^{t,y,z}=\lambda y-\int_{t}^{s}g(r,\lambda Y_{r}^{t,y,z}%
,Z_{r}^{t,y,z,\lambda})dr+\int_{t}^{s}Z_{r}^{t,y,z,\lambda}dB_{r}%
,\ s\in \lbrack t,T].
\]
Compare this with $\lambda Y_{s}^{t,y,z}=\lambda y-\int_{t}^{s}\lambda
g(r,Y_{r}^{t,y,z},z)dr+\int_{t}^{s}\lambda zdB_{r}$, it follows that
$Z_{\cdot}^{t,y,z,\lambda}\equiv \lambda z$ and $\lambda g(r,Y_{r}%
^{t,y,z},z)\equiv g(r,\lambda Y_{r}^{t,y,z},z)$, $r\in \lbrack t,T]$. We then
can apply Lemma \ref{Lem-f-f-} to obtain (i). $\blacksquare$

\medskip

\textbf{Proof of Proposition \ref{Proposition-4} }We first prove the
\textquotedblleft If\textquotedblright \ part. For each fixed $(y,z)\in R\times
R^{d}$, we have $Y_{s}^{t,y,z}\equiv \mathbb{E}_{s,T}^{g}[Y_{T}^{t,y,z}]\equiv
y+\mathbb{E}_{s,T}^{g}[Y_{T}^{t,y,z}-y]$. Let $\bar{Y}_{s}=\mathbb{E}%
_{s,T}^{g}[Y_{T}^{t,y,z}-y]$, $s\in \lbrack0,T]$ and $\bar{Z}$ be the
corresponding part of It\^{o}'s integrand. By $\bar{Y}_{r}\equiv
y+Y_{r}^{t,y,z}$ it follows that
\begin{align*}
y+Y_{s}  &  =y+Y_{T}^{t,y,z}+\int_{s}^{T}g(r,Y_{r}^{t,y,z},z)-\int_{s}%
^{T}zdB_{r}\\
&  =(y+Y_{T}^{t,y,z})+\int_{s}^{T}g(r,\bar{Y}_{r},\bar{Z}_{r})-\int_{s}%
^{T}\bar{Z}_{r}dB_{r}.
\end{align*}
Thus $\bar{Z}_{r}\equiv z$ and
\[
g(r,Y_{r}^{t,y,z},z)\equiv g(r,Y_{r}^{t,y,z}-y,\bar{Z}_{r})\equiv
g(r,Y_{r}^{t,y,z}-y,z).
\]
We then can apply Lemma \ref{Lem-f-f-} to obtain that, for each $(y,z)\in
R\times R^{d}$,
\[
g(r,y,z)\equiv g(r,y-y,z)\equiv g(r,0,z).
\]
Namely, $g$ is independent of $y$. \newline \textquotedblleft Only if
part\textquotedblright: For each for each $s\leq t$ and $X\in L^{2}%
(\mathcal{F}_{t})$, $\eta \in L^{2}(\mathcal{F}_{s})$, we have
\[
Y_{r}:=\mathbb{E}_{s,t}^{g}[X+\eta]=X+\eta+\int_{r}^{t}g(u,Z_{u})du-\int
_{s}^{t}Z_{u}dB_{u},\;r\in \lbrack s,t].
\]
Thus $\bar{Y}_{r}:=Y_{r}-\eta$ is a $g$--solution on $[s,t]$ with terminal
condition $\bar{Y}_{t}=X+\eta$. This implies
\[
\mathbb{E}_{s,t}^{g}[X]+\eta=\bar{Y}_{s}=\mathbb{E}_{s,t}^{g}[X+\eta].
\]
The proof is complete. $\blacksquare$

\medskip

\textbf{Proof of Proposition \ref{Proposition-5} }The \textquotedblleft
if\textquotedblright \ part: Since process $Y_{t}:=\mathbb{E}_{t,T}^{g}[X]$
solves the following BSDE
\[
Y_{t}=X+\int_{t}^{T}g(s,Y_{s},Z_{s})ds-\int_{t}^{T}Z_{s}dB_{s},
\]
we have
\[
Y_{t}+\int_{0}^{t}\bar{z}_{s}^{i_{0}}dB_{s}^{i_{0}}=X+\int_{0}^{T}\bar{z}%
_{s}^{i_{0}}dB_{s}^{i_{0}}+\int_{t}^{T}g(s,Y_{s},Z_{s})ds-\int_{t}^{T}\bar
{Z}_{s}dB_{s},
\]
where
\[
\bar{Z}_{s}=\left(  Z_{s}^{1},\cdots,Z_{s}^{i_{0}-1},Z_{s}^{i_{0}}+\bar{z}%
_{s}^{i_{0}},Z_{s}^{i_{0}+1},\cdots,Z_{s}^{d}\right)  .
\]
But since $g(s,y,z)$ does not depend the $i_{0}$--th component of $z\in R^{d}%
$, we thus have $g(s,Y_{s},Z_{s})\equiv g(s,Y_{s},\bar{Z}_{s})$. Thus
\[
Y_{t}+\int_{0}^{t}\bar{z}_{s}^{i_{0}}dB_{s}^{i_{0}}=X+\int_{0}^{t}\bar{z}%
_{s}^{i_{0}}dB_{s}^{i_{0}}+\int_{t}^{T}g(s,Y_{s},\bar{Z}_{s})ds-\int_{t}%
^{T}\bar{Z}_{s}dB_{s}.
\]
This means that (\ref{zBi0}) holds. \newline The \textquotedblleft only
if\textquotedblright \ part: For each fixed $(t,y,z)$, let $(Y_{s}%
^{t,y,z})_{s\geq t}$ be the solution of (\ref{Ytyz}) with $f=g$. We have,
\[
\mathbb{E}_{s,T}[Y_{T}^{t,y,z}]-z^{i_{0}}B_{s}^{i_{0}}=\mathbb{E}_{s,T}%
[Y_{T}^{t,y,z}-z^{i_{0}}B_{T}^{i_{0}}],\;s\in \lbrack t,T].
\]
Since the process $Y_{r}:=\mathbb{E}_{s,r}[Y_{r}^{t,y,z}-z_{r}^{i_{0}}%
B_{r}^{i_{0}}]$, $r\in \lbrack t,s]$, solves the BSDE
\[
Y_{s}^{t,y,z}-z^{i_{0}}B_{s}^{i_{0}}=Y_{s}=Y_{T}^{t,y,z}+z^{i_{0}}B_{T}%
^{i_{0}}+\int_{s}^{T}g(r,Y_{r},Z_{r})ds-\int_{s}^{T}Z_{r}dB_{r}.
\]
From which we deduce $Z_{s}=\bar{z}:=\left(  z^{1},\cdots,z^{i_{0}%
-1},0,z^{i_{0}+1},\cdots,z^{d}\right)  =z$ and thus
\[
g(r,Y_{r},Z_{r})=g(r,Y_{r}^{t,y,z},\bar{z})=g(r,Y_{r}^{t,y,z},z),\;0\leq t\leq
r\leq T.
\]
It then follows from Lemma \ref{Lem-f-f-} that
\[
g(t,y,\bar{z})=g(t,y,z),\;t\geq0\text{, a.e., a.s.,}%
\]
i.e., $g$ does not depend the $i_{0}$--th component of $z\in R^{d}$.
$\blacksquare$

\subsection{Testing the criteria (A5) by market data}

With Chen L. and Sun P. of our research group, we proceed a data test for the
criteria (A5), i.e., the domination inequality (\ref{domination}), to check if
a specific pricing mechanism is a $g$--expectation, or $g$--pricing mechanism
$\mathbb{E}^{g}$.

We have firstly tested the CME (Chicago Mercantile Exchange)'s market pricing
mechanism of derivatives by taking the daily closing prices of options with
S\&P500 index futures as the underlying asset. The data is obtained from
parameter files published from CME's fpt-webset, named cmeMMDDs.par (MM for
month, DD for day) of call and put prices, from\ 05 January 2000 to November
2003, of totally 960 trading days. The corresponding S\&P500 future's prices
is obtained from the parameter files of SPAN (Standard Portfolio Analysis of
Risk) system downloaded from CME's ftp site.

We denote by $X_{T}^{i}=(S_{T}-k_{i})^{+}$ (resp. $Y_{T}^{i}=(S_{T}-k_{i}%
)^{-}$), the market maturity value of the call (resp. put) option with
maturity $T$ and strike price $k_{i}$. The corresponding \ values of the short
positions are $-X_{T}^{i}$ and $-Y_{T}^{i}$. We denote the market price of the
corresponding prices of options at time $t<T$ by $\mathbb{E}_{t,T}^{m}%
[X_{T}^{i}]$, $\mathbb{E}_{t,T}^{m}[Y_{T}^{i}]$, $\mathbb{E}_{t,T}^{m}%
[-X_{T}^{i}]$ and $\mathbb{E}_{t,T}^{m}[-Y_{T}^{i}]$ respectively. The
inequalities we need to put to the test are, according to (\ref{domination}),
in the following different combinations, with different $(t,T)$ and different
strike prices%
\begin{equation}
\left \{
\begin{array}
[c]{ccc}%
\text{Call--Call:} &  & \mathbb{E}_{t,T}^{m}[X_{T}^{i}]-\mathbb{E}_{t,T}%
^{m}[X_{T}^{j}]\leq \mathbb{E}_{t,T}^{g_{\mu}}[X_{T}^{i}-X_{T}^{j}]\\
\text{Put--Put:} &  & \mathbb{E}_{t,T}^{m}[Y_{T}^{i}]-\mathbb{E}_{t,T}%
^{m}[-Y_{T}^{j}]\leq \mathbb{E}_{t,T}^{g_{\mu}}[Y_{T}^{i}-Y_{T}^{j}]\\
\text{Call--Put:} &  & \mathbb{E}_{t,T}^{m}[X_{T}^{i}]-\mathbb{E}_{t,T}%
^{m}[Y_{T}^{j}]\leq \mathbb{E}_{t,T}^{g_{\mu}}[X_{T}^{i}-Y_{T}^{j}]\\
\text{Put--Call:} &  & \mathbb{E}_{t,T}^{m}[Y_{T}^{i}]-\mathbb{E}_{t,T}%
^{m}[X_{T}^{j}]\leq \mathbb{E}_{t,T}^{g_{\mu}}[Y_{T}^{i}-X_{T}^{j}]
\end{array}
\right.  \label{testedIQ}%
\end{equation}
and%
\begin{equation}
\left \{
\begin{array}
[c]{ccc}%
\text{Call--ShortCall:} &  & \mathbb{E}_{t,T}^{m}[X_{T}^{i}]-\mathbb{E}%
_{t,T}^{m}[-X_{T}^{j}]\leq \mathbb{E}_{t,T}^{g_{\mu}}[X_{T}^{i}+X_{T}^{j}]\\
\text{Put--ShortPut:} &  & \mathbb{E}_{t,T}^{m}[Y_{T}^{i}]-\mathbb{E}%
_{t,T}^{m}[-Y_{T}^{j}]\leq \mathbb{E}_{t,T}^{g_{\mu}}[Y_{T}^{i}+Y_{T}^{j}]\\
\text{Call--ShortPut:} &  & \mathbb{E}_{t,T}^{m}[X_{T}^{i}]-\mathbb{E}%
_{t,T}^{m}[-Y_{T}^{j}]\leq \mathbb{E}_{t,T}^{g_{\mu}}[X_{T}^{i}+Y_{T}^{j}]
\end{array}
\right.  \label{testedIQs}%
\end{equation}
\  \ In the above inequalities the data of the left hand sides is the market
prices of options taken from CME parameter files. In our testing the
transaction cost is neglected, i.e., we assume that $\mathbb{E}_{t,T}%
^{m}[-X]=-\mathbb{E}_{t,T}^{m}[X]$. The right hand sides is the corresponding
values of $g_{\mu}$--expectations. We fix $\mu=25$ uniformly for all tested
inequalities. We have calculated all these values on the right hand side by
using standard binomial tree algorithm of BSDE. Here an improved version of
the algorithms of BSDE proposed Peng and Xu [2005] has been applied to solve
the following 1-dimensional BSDE:%

\begin{align}
y_{t} &  =X+\int_{t}^{T}\mu(|y_{s}|+|z_{s}|)ds-\int_{t}^{T}z_{s}%
dB_{s}\label{Numerical}\\
&  \text{\ }\nonumber
\end{align}
with different terminal values $y_{T}=X_{T}^{i}-X_{T}^{j}$, $Y_{T}^{i}%
-Y_{T}^{j}$,\ $X_{T}^{i}-Y_{T}^{j}$, $Y_{T}^{i}-X_{T}^{j}$, $X_{T}^{i}%
+X_{T}^{j}$, $Y_{T}^{i}+Y_{T}^{j}$,\ $X_{T}^{i}+Y_{T}^{j}$, respectively. The
closing prices of S\&P500 futures options of 69 trading days from year 2000 to
2003 have been put in the test. With the above mentioned combinations, we have
tested a total number of 6,200,828 inequalities of (\ref{testedIQ}) and
(\ref{testedIQs}). \ This means that our BSDE (\ref{Numerical}) have been
calculated 6,200,828 times. A very positive result was obtained: among the
totally 6,200,828 tested inequalities, only 17 are against the criteria
(\ref{testedIQ}). Among those 12 cases of violations, 5 are singular situation
since they themselves all violate Axiomatic monotonicity condition (A1). 5
cases are all from the same file cme0701s.par, 2003, Put--Put. They are all
the following singular cases:
\[
\mathbb{E}_{t,T}^{m}[(S_{T}-k_{i})^{-}]>\mathbb{E}_{t,T}^{m}[(S_{T}-k_{j}%
)^{-}],\  \  \text{for }k_{i}>k_{j}.
\]
The other 12 violations are the cases where the time $T-t$ is too short (less
than 2 days).

Since we have not found available data of bid-ask prices of the above options
from CME, we then have tested the bid-ask pricing mechanism of S\&P500 index
options operated by the system of market makers of CBOE The data source is
from Yahoo's finance quotes of the option prices from 07 December to 08 May
2006. We have collected the prices of 5,000 time points, i.e., 5,000 different
$t$ among 100 trading days. We denote this pricing mechanism by $\mathbb{E}%
_{t,T}^{mm}[X]$ for the ask price of an option $X$. According to our point of
view the bid price of the same $X$ is $-\mathbb{E}_{t,T}^{mm}[-X]$ and thus
the bid--ask spread is $\mathbb{E}_{t,T}^{mm}[X]+\mathbb{E}_{t,T}^{mm}[-X]$.
We have tested a total number of 589,360 inequalities of (\ref{testedIQ}) and
(\ref{testedIQs}), with $\mathbb{E}^{mm}$ in the place of $\mathbb{E}^{m}$.
Only $1$ case of violation appears.

We will report these test results in details in our forthcoming paper.
\cite{ChenSunPeng}.

\end{document}